\documentclass{amsart}
\usepackage{amsmath,amsfonts,latexsym,amssymb}

\newtheorem{theorem}{Theorem}
\newtheorem{lemma}{Lemma}
\newtheorem{corollary}{Corollary}
\newtheorem{proposition}{Proposition}

\theoremstyle{remark}
\newtheorem{remark}{Remark}

\begin{document}

\title[Determination of $GL(3)$ Hecke-Maass forms]{Determination of $GL(3)$ Hecke-Maass forms\\ from twisted central values}
\author{Ritabrata Munshi \and Jyoti Sengupta}
\address{School of Mathematics, Tata institute of Fundamental Research, 1 Dr. Homi Bhabha Road, Colaba, Mumbai 400005, India.} 
\email{rmunshi@math.tifr.res.in \and sengupta@math.tifr.res.in}


\begin{abstract}
Suppose $\pi_1$ and $\pi_2$ are two Hecke-Maass cusp forms for $SL(3,\mathbb{Z})$ such that for all primitive character $\chi$ we have 
$$
L(\tfrac{1}{2},\pi_1\otimes\chi)=L(\tfrac{1}{2},\pi_2\otimes\chi).
$$
Then we show that $\pi_1=\pi_2$. 
\end{abstract}

\maketitle

%======================================================================================================================
%======================================================================================================================

\section{Introduction}
\label{intro}

Determining modular forms  from central values of the $L$ 
function of its twists has a fairly long history (see \cite{CD}, \cite{GHS}, \cite{Li}, \cite{Liu}, \cite{Lu}, \cite{LR}, \cite{M}, \cite{MS}, \cite{P} and \cite{Z}) and it still remains a topic of much interest.  It was first considered by Luo and Ramakrishnan \cite{LR}.  They showed that if two cuspidal normalised newforms $f$ and $g$ of weight $2k$ (resp. $2k^{'}$) and level $N$
(resp. $N^{'}$) have the property that
$$
L (\tfrac{1}{2}, f \otimes \chi) = L(\tfrac{1}{2}, g \otimes \chi)
$$
for all characters $\chi$, then $k=k^{'}, \ N=N^{'}$ and $f=g$. (In the same paper they prove a much stronger result where only quadratic twists are required.) The aim of this article is to generalise this result to $GL(3)$ Maass forms of full level.  More precisely we prove the following.\\

\begin{theorem}
\label{mthm}
Let $\pi_1$ and $\pi_2$ be two Hecke-Maass cusp forms for $SL(3,\mathbb{Z})$. Suppose for all primitive character $\chi$ we have the equality of the central values of the twisted $L$-functions
$$
L(\tfrac{1}{2},\pi_1\otimes\chi)=L(\tfrac{1}{2},\pi_2\otimes\chi).
$$
Then $\pi_1=\pi_2$.
\end{theorem}

\begin{remark}
When $\pi_1$ and $\pi_2$ are self dual and $\chi$ is 
restricted  to quadratic characters $\chi_d$, the above result has been 
established by Chinta and Diaconu.  They use the method of multiple Dirichlet
series which is completely different from our method which employs twisted
averages.
\end{remark}

As in the case of $GL(2)$ forms the result will follow from suitable asymptotic for the twisted central values. But the task of computing the first moment of $GL(3)$ $L$-functions over a family of twists is much more delicate. In fact we can successfully do it only in the special case of family of twists by characters of almost prime modulus. More precisely, we take 
$$
\mathcal{Q}=\{q=q_1q_2:\; Q_1<q_1<2Q_1,\;Q_2<q_2<2Q_2,\;q_1, q_2\;\text{primes}\}
$$
with $Q_1=Q^{3/4-\delta}$, $Q_2=Q^{1/4+\delta}$, and $\delta=1/100$. Our main task in this paper is to prove the following proposition.\\

\begin{proposition}
\label{prop}
For $\mathcal{Q}$ as above, we have
\begin{align*}
\sum_{q\in\mathcal{Q}}\;\sideset{}{^\star}\sum_{\chi\bmod{q}}L\left(\tfrac{1}{2},\pi\otimes\chi\right)(1+\chi(-1))\bar\chi(\ell)=\frac{\lambda_\pi(\ell,1)}{\sqrt{\ell}}Y+O\left(\ell^2\:Q^{2-1/2013+\varepsilon}\right),
\end{align*}
where $Y=\sum_{q\in\mathcal{Q}}q$ and the implied constant is independent of $\pi$.
Here the $\star$ indicates that we are summing only over primitive characters $\chi$.
\end{proposition}

A similar idea of averaging over factorizable moduli was also utilized by Luo \cite{Luo}. As in \cite{Luo} we also use unbalanced approximate functional equation and Deligne's bound for hyper-Kloosterman sums. There are however some subtle differences between this work and \cite{Luo}, and the Deligne's bound for some more complicated exponential sums are required in this work.\\

We conclude this section by deriving Theorem~\ref{mthm} from Proposition~\ref{prop}. Applying the proposition we see that the equality of the central values of the twisted $L$-functions implies equality among the Fourier coefficients
$$
\lambda_{\pi_1}(\ell,1)=\lambda_{\pi_2}(\ell,1).
$$
The strong multiplicity one theorem of Jacquet and Shalika (see \S 12.6 of \cite{G}) then implies that $\pi_1=\pi_2$. The theorem follows.\\

%===============================================================================
\section{The set up}
Let $\pi$ be a $SL(3,\mathbb{Z})$ Hecke-Maass cusp form with (Hecke) normalized Fourier coefficients $\lambda_\pi(n,m)$. Let $\chi$ be an even ($\chi(-1)=1$) primitive character modulo $q$. The $L$-function associated with the twisted form $\pi\otimes\chi$ is given by the Dirichlet series
$$
L(s,\pi\otimes\chi)=\sum_{n=1}^\infty\frac{\lambda_\pi(n,1)\chi(n)}{n^s}
$$
in the half plane $\text{Re}(s)=\sigma>1$. The $L$-function extends to an entire function and satisfies the functional equation
\begin{align*}
\Lambda(s,\pi\otimes\chi)=\varepsilon_\chi \Lambda(1-s,\tilde{\pi}\otimes\bar{\chi})
\end{align*}
where 
\begin{align*}
\Lambda(s,\pi\otimes\chi)=q^{3s/2}\gamma(s)L(s,\pi\otimes\chi)
\end{align*}
is the completed $L$-function. The gamma factor $\gamma(s)$ can be explicitly expressed in terms of the Langlands parameters of $\pi$, namely 
$$
\gamma(s)=\prod_{i=1}^3\Gamma_{\mathbb{R}}(s-\alpha_i)
$$
with $\Gamma_{\mathbb{R}}(s)=\pi^{-s/2}\Gamma(s/2)$. (Since we are restricting ourselves to even characters $\gamma(s)$ is independent of $\chi$. In general $\gamma(s)$ also depends on the parity of $\chi$.) The explicit expression of $\gamma$ will not be required in this work. But we will use that fact that $\text{Re}(\alpha_i)\leq \frac{2}{5}$ (see \cite{LRS}). The sign of the functional equation is given by
\begin{align*}
\varepsilon_\chi=g(\chi)^3/q^{3/2}
\end{align*}
where $g(\chi)$ is the associated Gauss sum. Note that $\tilde{\pi}$ is the dual form, and $\bar{\chi}$ denotes the complex conjugate. Also note that 
$$
\lambda_{\tilde{\pi}}(n,m)=\overline{\lambda_\pi(n,m)}=\lambda_\pi(m,n).
$$\\

The approximate functional equation gives
\begin{align}
\label{afe}
L\left(\tfrac{1}{2},\pi\otimes\chi\right)=\sum_{n=1}^\infty \frac{\lambda_\pi(n,1)\chi(n)}{\sqrt{n}}V\left(\frac{nX}{q^{\frac{3}{2}}}\right)+\frac{g(\chi)^3}{q^{\frac{3}{2}}}\sum_{n=1}^\infty \frac{\lambda_{\tilde\pi}(n,1)\bar\chi(n)}{\sqrt{n}}V\left(\frac{n}{Xq^{\frac{3}{2}}}\right)
\end{align}
where the function $V$ is given by  
$$
V(y)=\frac{1}{2\pi i}\int_{(3)}y^{-s}\frac{\gamma(s+1/2)}{\gamma(1/2)}\frac{\mathrm{d}s}{s}.
$$
The parameter $X>0$ will be optimally chosen later. For $y$ small we can shift the contour to the left upto $\sigma=-1/10+\varepsilon$, to get the asymptotic expansion
\begin{align}
\label{exp-v}
V(y)=1+O(y^{1/10-\varepsilon}).
\end{align} 
The implied constant depends only on $\varepsilon$.\\

We want to compute the twisted average
$$
\mathfrak{T}=\sum_{q\in\mathcal{Q}}\;\sideset{}{^\star}\sum_{\chi\bmod{q}}L\left(\tfrac{1}{2},\pi\otimes\chi\right)(1+\chi(-1))\bar\chi(\ell),
$$
where $\chi$ ranges over the collection of primitive (even) characters modulo $q$. We take $\ell$ to be a small power of a prime, and $\mathcal{Q}$ is a collection of square-free integers in the interval $[Q,2Q]$. Using \eqref{afe} we get
$$
\mathfrak{T}=\mathfrak{F}+\mathfrak{S}
$$
where
\begin{align}
\label{f}
\mathfrak{F}=\sum_{q\in\mathcal{Q}}\sum_{n=1}^\infty \frac{\lambda_\pi(n,1)}{\sqrt{n}}V\left(\frac{nX}{q^{\frac{3}{2}}}\right)\;\sideset{}{^\star}\sum_{\chi\bmod{q}}\chi(n)(1+\chi(-1))\bar\chi(\ell)
\end{align}
and
\begin{align}
\label{s}
\mathfrak{S}=\sum_{q\in\mathcal{Q}}\frac{1}{q^{\frac{3}{2}}}\sum_{n=1}^\infty \frac{\lambda_{\tilde\pi}(n,1)}{\sqrt{n}}V\left(\frac{n}{Xq^{\frac{3}{2}}}\right)\sideset{}{^\star}\sum_{\chi\bmod{q}}g(\chi)^3(1+\chi(-1))\bar\chi(n\ell).
\end{align}\\

We will consider the case 
$$
\mathcal {Q}=\mathcal{Q}_1\mathcal{Q}_2,
$$ 
where $\mathcal{Q}_i$ is the collection of prime numbers $q_i\equiv 1\bmod{4}$ in the range $[Q_i,2Q_i]$. We take $2Q_1<Q_2$. So each $q\in\mathcal{Q}$ can be uniquely written as $q_1q_2$ where $q_i\in\mathcal Q_i$. 
Note that a character $\chi$ modulo $q=q_1q_2$ is primitive if and only if it splits as $\chi=\chi_1\chi_2$ with $\chi_i$ primitive modulo $q_i$.\\

%==========================================================================
\section{Character sums}
Since $q_i$ is a prime number, the only non-primitive character modulo $q_i$ is the trivial character. So it follows that
\begin{align}
\label{identity-1}
\sideset{}{^\star}\sum_{\chi_i\bmod{q_i}}\chi_i(n)\bar{\chi_i}(\ell)=\phi(q_i)1_{n\equiv \ell\bmod{q_i}}-1.
\end{align}
This identity will be required to evaluate the contribution of the first term of the approximate functional equation \eqref{afe}. To evaluate the contribution of the second term we need to tackle a more complicated sum because of the involvement of the Gauss sums coming from the sign of the functional equation.\\

We first consider the Gauss sum of the product $\chi_1\chi_2$, i.e.
$$
g(\chi_1\chi_2)=\sum_{a\bmod{q_1q_2}}\chi_1(a)\chi_2(a)e_{q_1q_2}(a).
$$ 
Each $a$ in the above sum can be written uniquely as 
$$
a=a_1q_2\bar{q}_2+a_2q_1\bar{q}_1
$$
with $a_i\bmod{q_i}$. Consequently we get
\begin{align*}
g(\chi_1\chi_2)&=\sum_{a_1\bmod{q_1}}\sum_{a_2\bmod{q_2}}\chi_1(a_1)\chi_2(a_2)e_{q_1}(a_1\bar{q}_2)e_{q_2}(a_2\bar{q}_1)\\
&=\chi_1(q_2)\chi_2(q_1)g(\chi_1)g(\chi_2).
\end{align*}
This leads us to consider the sum
\begin{align*}
\sideset{}{^\star}\sum_{\chi_i\bmod{q_i}}g(\chi_i)^3\chi_i(r)\bar{\chi_i}(n\ell),
\end{align*}
with $(r,q_i)=1$.
Opening the Gauss sum we get
\begin{align*}
\sideset{}{^\star}\sum_{\chi_i\bmod{q_i}}\left(\sum_{a\bmod{q_i}}\chi_i(a)e_{q_i}(a)\right)^3\chi_i(r)\bar{\chi_i}(n\ell).
\end{align*}
Then interchanging the order of summations and using \eqref{identity-1} we arrive at
\begin{align*}
&\sideset{}{^\star}\sum_{a,b,c\bmod{q_i}}e_{q_i}(a+b+c)\sideset{}{^\star}\sum_{\chi_i\bmod{q_i}}\chi_i(rabc)\bar{\chi_i}(n\ell)\\
=&\sideset{}{^\star}\sum_{a,b,c\bmod{q_i}}e_{q_i}(a+b+c)\left\{\phi(q_i)1_{rabc\equiv n\ell\bmod{q_i}}-1\right\}\\
=&\phi(q_i)\sideset{}{^\star}\sum_{a,b\bmod{q_i}}e_{q_i}(a+b+n\ell\overline{rab})-c_{q_i}(1)^3.
\end{align*}
Here $c_{q_i}(1)$ stands for the Ramanujan sum modulo $q_i$, and since we are taking $q_i$ prime we have $c_{q_i}(1)=-1$. We set
\begin{align}
\label{kloosterman}
K_{q_i}(u)=\sideset{}{^\star}\sum_{a,b\bmod{q_i}}e_{q_i}(a+b+u\overline{ab}).
\end{align}
This is a hyper-Kloosterman sum. Square-root cancellation in such sums was conjectured by Mordell. Deligne later proved that 
$$
K_{q_i}(u)\ll q_i
$$
for any integer $u$. Also note that 
\begin{align}
\label{inv}
\overline{K_{q_i}(u)}=K_{q_i}(-u).
\end{align} 
We conclude the following lemma.\\

\begin{lemma}
We have
\begin{align}
\label{identity-2}
\sideset{}{^\star}\sum_{\chi_i\bmod{q_i}}g(\chi_i)^3\chi_i(r)\bar{\chi_i}(n\ell)=\phi(q_i)K_{q_i}(nl\bar{r})+1.
\end{align}
\end{lemma}

\bigskip

%=============================================================================

\section{The first term}
In this section we will establish an asymptotic expression for the first term \eqref{f}. First using \eqref{identity-1} we get
\begin{align*}
\sideset{}{^\star}\sum_{\chi\bmod{q}}\chi(n)(1+\chi(-1))\bar\chi(\ell)&=\sum_\pm\sideset{}{^\star}\sum_{\chi\bmod{q}}\chi(n)\bar\chi(\pm\ell)\\
=\sum_{\pm}\prod_{i=1}^2\left(\sideset{}{^\star}\sum_{\chi_i\bmod{q_i}}\chi_i(n)\bar\chi_i(\pm\ell)\right)&=
\sum_{\pm}\prod_{i=1}^2\left(\phi(q_i)1_{n\equiv \pm\ell\bmod{q_i}}-1\right).
\end{align*}
Accordingly $\mathfrak{F}$ splits as a sum of eight terms
$$
\mathfrak{F}=\sum_\pm\sum_{r\in R}\:(-1)^{\nu(r)}\:\mathfrak{F}_{\pm,r}
$$
where $R=\{1,q_1,q_2,q_1q_2\}$, and
\begin{align*}
\mathfrak{F}_{\pm,r}=\sum_{q\in\mathcal{Q}}\;\phi(r)\sum_{n\equiv \ell\bmod{r}} \frac{\lambda_\pi(n,1)}{\sqrt{n}}V\left(\frac{nX}{q^{\frac{3}{2}}}\right).
\end{align*}
Here $\nu(r)$ denotes the number of prime factors of $r$.
Recall that we have $\mathcal{Q}=\mathcal{Q}_1\mathcal{Q}_2$ and so $q=q_1q_2$, with $q_i$ primes.\\

%Let us consider the case where $r=1$. To this term we apply the Cauchy inequality and the bound
%$$
%\sum_{n\leq N}|\lambda_\pi(n,1)|^2\ll N^{1+\varepsilon}
%$$ 
%to get
%$$
%\mathfrak{F}_{\pm,1}\ll \sum_{q\in\mathcal{Q}}\;\sum_{n=1}^\infty \frac{|\lambda_\pi(n,1)|}{\sqrt{n}}V\left(\frac{nX}{q^{\frac{3}{2}}}\right)\ll Q^{7/4+\varepsilon}X^{-1/2}.
%$$\\

We will evaluate
$$
\mathfrak{F}_{\pm,r}=\sum_{q\in\mathcal{Q}}\;\phi(r)\sum_{n\equiv \pm\ell\bmod{r}} \frac{\lambda_\pi(n,1)}{\sqrt{n}}V\left(\frac{nX}{q^{\frac{3}{2}}}\right),
$$
by changing the order of summations 
$$
\mathfrak{F}_{\pm,r}=\sum_{n=1}^\infty \frac{\lambda_\pi(n,1)}{\sqrt{n}}V\left(\frac{nX}{q^{\frac{3}{2}}}\right)\sum_{\substack{q\in\mathcal{Q}\\r|n\mp \ell}}\;\phi(r).
$$
For any given $n$, the sum over $q$ is restricted by the condition $r|n\mp\ell$. Naturally the diagonal $n=\ell$ yields the leading term for $\mathfrak{F}_{+,r}$, namely
$$
\mathfrak{L}_r=\frac{\lambda_\pi(\ell,1)}{\sqrt{\ell}}V\left(\frac{\ell X}{q^{\frac{3}{2}}}\right)\:\sum_{q\in\mathcal{Q}}\;\phi(r).
$$
Since $\ell>0$, there are no leading terms for $\mathfrak{F}_{-,r}$. Now
$$
\sum_{r\in R}\mathfrak{L}_r=\frac{\lambda_\pi(\ell,1)}{\sqrt{\ell}}V\left(\frac{\ell X}{q^{\frac{3}{2}}}\right)\:\sum_{q\in\mathcal{Q}}\;q.
$$
Set
$$
Y=\sum_{q\in\mathcal{Q}}\;q
$$
which is roughly of size $Q^2/\log Q_1\log Q_2$. We want a main term which will involve only the Fourier coefficient of the form $\pi$, and no other parameters. The above sum still involves the Langlanda parameters of $\pi$. We now use \eqref{exp-v} to arrive at
$$
\sum_{r\in R}\mathfrak{L}_r=\frac{\lambda_\pi(\ell,1)}{\sqrt{\ell}}Y+O\left(\ell^{3/5}X^{1/10}\:Q^{2-3/20+\varepsilon}\right).
$$
Here we have used the trivial bound $|\lambda(\ell,1)|\ll \ell$.\\

Next we study the off-diagonal $n\neq \pm\ell$. In this case
$$
\sum_{\substack{q\in\mathcal{Q}\\r|n\mp \ell}}\;\phi(r)\ll Q^{1+\varepsilon}.
$$
Hence the contribution of this term to $\mathfrak{F}_{\pm,r}$ is dominated by
\begin{align}
\label{as-in}
Q^{1+\varepsilon}\sum_{n=1}^\infty \frac{|\lambda_\pi(n,1)|}{\sqrt{n}}\Bigr|V\left(\frac{nX}{q^{\frac{3}{2}}}\right)\Bigl|.
\end{align}
We apply the Cauchy inequality and the well-known bound
\begin{align}
\label{ram-on-av}
\sum_{n\leq N}|\lambda_\pi(n,1)|^2\ll N^{1+\varepsilon}
\end{align}
(which follows from the Rankin-Selberg theory) to get that the above sum is dominated by
$$
O\left(Q^{7/4+\varepsilon}/\sqrt{X}\right).
$$
\\

\begin{lemma}
\label{lem-f}
We have
$$
\mathfrak{F}=\frac{\lambda_\pi(\ell,1)}{\sqrt{\ell}}Y+O\left(\ell^{3/5}X^{1/10}\:Q^{2-3/20+\varepsilon}+
Q^{7/4+\varepsilon}/\sqrt{X}\right)
$$
where $Y=\sum_{q\in\mathcal{Q}}q$.
\end{lemma}

We note that the last error term is satisfactory for our purpose as long as $X\gg Q^{-1/2+\delta}$ for any $\delta>0$. Indeed by choosing $X=Q^{-1/2+\delta}$ the above error term reduces to
\begin{align}
\label{err-term-1}
O\left(\ell^{3/5}\:Q^{9/5+\delta/10+\varepsilon}+
Q^{2-\delta/2+\varepsilon}\right).
\end{align}\\

%===================================================================================
\section{The second term}
Next we study the contribution of the second term \eqref{s}. Using \eqref{identity-2} we get
\begin{align*}
&\sideset{}{^\star}\sum_{\chi\bmod{q}}g(\chi)^3(1+\chi(-1))\bar\chi(n\ell)=\sum_\pm \sideset{}{^\star}\sum_{\chi\bmod{q}}g(\chi)^3\bar\chi(\pm n\ell)\\
&=\sum_\pm \left(\sideset{}{^\star}\sum_{\chi_1\bmod{q_1}}g(\chi_1)^3\chi_1(q_2^3)\bar\chi_1(\pm n\ell)\right)\left(\sideset{}{^\star}\sum_{\chi_2\bmod{q_2}}g(\chi_2)^3\chi_2(q_1^3)\bar\chi_2(\pm n\ell)\right)\\
&=\sum_\pm\left(\phi(q_1)K_{q_1}(\pm nl\bar{q}_2^3)+1\right)\left(\phi(q_2)K_{q_2}(\pm nl\bar{q}_1^3)+1\right).
\end{align*}
Now by Deligne's bound for hyper-Kloosterman sums the above expression reduces to
$$
\sum_\pm\:\phi(q)K_{q_1}(\pm n\ell\bar{q}_2^3)K_{q_2}(\pm n\ell\bar{q}_1^3)+O(Q(Q_1+Q_2)).
$$\\

The total contribution of the error term to \eqref{s} is bounded by
\begin{align*}
\frac{(Q_1+Q_2)}{\sqrt{Q}}\:\sum_{q\in\mathcal{Q}}\;\sum_{n=1}^\infty \frac{|\lambda_{\tilde \pi}(n,1)|}{\sqrt{n}}V\left(\frac{n}{q^{\frac{3}{2}}X}\right).
\end{align*}
Compare this with \eqref{as-in}. We conclude that this is dominated by
$$
O\left((Q_1+Q_2)Q^{5/4+\varepsilon}X^{1/2}\right).
$$
This term is satisfactory for our purpose if $\sqrt{X}(Q_1+Q_2)\ll Q^{3/4-\delta}$ for some $\delta>0$.\\

Finally we need to analyse the sum
$$
\mathfrak{R}_\pm=\sum_{q\in\mathcal{Q}}\;\frac{\phi(q)}{q^{3/2}}\sum_{n=1}^\infty \frac{\lambda_{\tilde \pi}(n,1)}{\sqrt{n}}K_{q_1}(\pm n\ell\bar{q}_2^3)K_{q_2}(\pm n\ell\bar{q}_1^3)V\left(\frac{n}{q^{\frac{3}{2}}X}\right).
$$
Using the Deligne bound for hyper-Kloosterman sums and the Ramanujan bound on average, we see that the above sum is bounded by 
$$
O(Q^{9/4+\varepsilon}\sqrt{X}).
$$
Since to control the error term from the first part of the functional equation we need to take $X>Q^{-1/2+\delta}$ with $\delta>0$, we see that this bound is worse than the main term which is of size $Q^2$ (roughly speaking). So we need non-trivial saving in the sum over $n$ and $q_i$. We can truncate the sum over $n$ at $Q^{3/2+\varepsilon}\sqrt{X}$ at a cost of a negligible error term. This follows from the rapid decay of the function $V$. In the remaining sum we substitute the explicit expression for $V$ and shift the contour to $\sigma=\varepsilon$ to get
$$
\mathfrak{R}_\pm=\frac{1}{2\pi i}\int_{(\varepsilon)}\frac{\gamma(s+1/2)}{\gamma(1/2)}X^{-s}\mathfrak{R}_\pm(s)\frac{\mathrm{d}s}{s}+O(Q^{-2013}),
$$
where
$$
\mathfrak{R}_\pm(s)=\sum_{q\in\mathcal{Q}}\;\frac{\phi(q)}{q^{3/2-s}}\sum_{1\leq n\ll N} \frac{\lambda_{\tilde \pi}(n,1)}{n^{1/2+s}}K_{q_1}(\pm n\ell\bar{q}_2^3)K_{q_2}(\pm n\ell\bar{q}_1^3),
$$
with $N=Q^{3/2+\varepsilon}\sqrt{X}$. Due to the rapid decay of the gamma factor $\gamma(s)$ we only need to get sufficient bound for $\mathfrak{R}_\pm(s)$ for $s=\varepsilon+it$ with $|t|\ll Q^\varepsilon$.\\

For notational simplicity we will focus on the $+$ term. The other term can be analysed in the same fashion. Taking absolute values we get
$$
\mathfrak{R}_+(s)\ll \frac{Q^\varepsilon}{\sqrt{Q_1}}\sum_{q_1\in\mathcal{Q}_1}\;\sum_{n\ll N} \frac{|\lambda_{\tilde \pi}(n,1)|}{\sqrt{n}}\Bigl|\sum_{q_2\in\mathcal{Q}_2}\;\frac{\phi(q_2)}{q_2^{3/2-s}}K_{q_1}( n\ell\bar{q}_2^3)K_{q_2}(n\ell\bar{q}_1^3)\Bigr|,
$$
for $s=\varepsilon+it$.
Next we apply the Cauchy inequality. Observe that from \eqref{ram-on-av}, it follows that 
$$
\sum_{q_1\in\mathcal{Q}_1}\;\sum_{n\ll Q^{3/2+\varepsilon}\sqrt{X}} \frac{|\lambda_{\tilde \pi}(n,1)|^2}{n}\ll Q_1Q^\varepsilon.
$$
Consequently
$$
\mathfrak{R}_+(s)\ll Q^\varepsilon\sqrt{\mathfrak{E}(s)}
$$
where
\begin{align}
\label{e}
\mathfrak{E}(s)=\sum_{q_1\in\mathcal{Q}_1}\;\sum_{1\leq n\ll N}\Bigl|\sum_{q_2\in\mathcal{Q}_2}\;\frac{\phi(q_2)}{q_2^{3/2-s}} K_{q_1}(n\ell\bar{q}_2^3)K_{q_2}(n\ell\bar{q}_1^3)\Bigr|^2.
\end{align}
Now observe that the same bound also holds for $\mathfrak{R}_-(\bar{s})$ because of \eqref{inv}.\\

\begin{lemma}
\label{lem-s}
We have
$$
\mathfrak{S}\ll Q^\varepsilon\int_{-Q^\varepsilon}^{Q^\varepsilon}\:\sqrt{\mathfrak{E}(\varepsilon+it)}\mathrm{d}t+(Q_1+Q_2)Q^{5/4+\varepsilon}X^{1/2}
$$
where $\mathfrak{E}(s)$ is as given in \eqref{e}.
\end{lemma}

\bigskip
%=================================================================================
\section{The remainder term}
We will now open the absolute square interchange the order of summations and apply Poisson summation on the sum over $n$ in \eqref{e}. But first we should smooth out the sum over $n$ using a smooth bump function $W$. This is possible as we have positivity. Indeed we have
\begin{align*}
\mathfrak{E}(s)\leq \sum_{q_1\in\mathcal{Q}_1}\;\sum_{n\in\mathbb{Z}}W\left(\frac{n}{N}\right)\Bigl|\sum_{q_2\in\mathcal{Q}_2}\;\frac{\phi(q_2)}{q_2^{3/2-s}} K_{q_1}(n\ell\bar{q}_2^3)K_{q_2}(n\ell\bar{q}_1^3)\Bigr|^2,
\end{align*}
for some compactly supported smooth function $W:\mathbb{R}\rightarrow\mathbb{R}$ satisfying $W^{(j)}(x)\ll_j 1$. We now conclude that
\begin{align}
\label{to-replace-1}
\mathfrak{E}(s)\ll \frac{Q^\varepsilon}{Q_2}\sum_{q_1\in\mathcal{Q}_1}\sum_{q_2\in\mathcal{Q}_2}\;\sum_{q_2'\in\mathcal{Q}_2}\;|\mathfrak{E}_{q_1,q_2,q_2'}|
\end{align}
where
$$
\mathfrak{E}_{q_1,q_2,q_2'}=\sum_{n\in\mathbb{Z}} K_{q_1}(n\ell\bar{q}_2^3)K_{q_2}(n\ell\bar{q}_1^3)K_{q_1}(-n\ell\bar{q_2'}^3)K_{q_2'}(-n\ell\bar{q}_1^3)W\left(\frac{n}{N}\right).
$$ \\

We split the sum over $n$ into congruence classes modulo $q_1q_2q_2'$ to get
\begin{align*}
\mathfrak{E}_{q_1,q_2,q_2'}=\sum_{\alpha\bmod{q_1q_2q_2'}}K_{q_1}(\alpha\ell\bar{q}_2^3)K_{q_2}(\alpha\ell\bar{q}_1^3)&K_{q_1}(-\alpha\ell\bar{q_2'}^3)K_{q_2'}(-\alpha\ell\bar{q}_1^3)\\
\times &\sum_{n\in \mathbb{Z}}W\left(\frac{\alpha+nq_1q_2q_2'}{N}\right).
\end{align*}
Apply the Poisson summation formula and make the change of variables $y=(\alpha+xq_1q_2q_2')/N$ to arrive at
\begin{align*}
\frac{N}{q_1q_2q_2'}\sum_{n\in\mathbb Z}&\mathfrak{C}(n;q_1,q_2,q_2')\int_{\mathbb R}W\left(y\right)e\left(-\frac{nNy}{q_1q_2q_2'}\right)\mathrm{d}y
\end{align*}
where the character sum is given by
\begin{align*}
\mathfrak{C}(n;q_1,q_2,q_2')=\sum_{\alpha\bmod{q_1q_2q_2'}}K_{q_1}(\alpha\ell\bar{q}_2^3)K_{q_2}(\alpha\ell\bar{q}_1^3)K_{q_1}(-\alpha\ell\bar{q_2'}^3)K_{q_2'}(-\alpha\ell\bar{q}_1^3)e_{q_1q_2q_2'}(\alpha n)
\end{align*}
Integrating by parts we see that the integral is negligibly small if $|n|\gg Q^{1+\varepsilon}Q_2/N$. Consequently
\begin{align}
\label{to-replace}
\mathfrak{E}_{q_1,q_2,q_2'}\ll \frac{\sqrt{QX}}{Q_2}\sum_{|n|\ll Q_2Q^\varepsilon/\sqrt{QX}}|\mathfrak{C}(n;q_1,q_2,q_2')|+Q^{-2013}.
\end{align}\\

%By our choice of $X$ we shall have $QX>Q^\delta$ for some $\delta>0$. So in the above sum, except for $n=0$ (the zero frequency), the other values of $n$ are necessarily co-prime with $q_2q_2'$. 
Consider the above character sum. Since $(q_1,q_2q_2')$, every $\alpha$ can be uniquely expressed as
$$
\alpha=\alpha_1 q_2q_2'\bar{q}_2\bar{q}_2'+\alpha_2q_1\bar{q}_1
$$
with $\alpha_1$ ranging modulo $q_1$ and $\alpha_2$ ranging modulo $q_2q_2'$. The character sum splits as a product of
\begin{align}
\label{a}
\mathfrak{A}(n;q_1,q_2,q_2')=\sum_{\alpha_1\bmod{q_1}}K_{q_1}(\alpha_1\ell\bar{q}_2^3)K_{q_1}(-\alpha_1\ell\bar{q_2'}^3)e_{q_1}(\alpha_1 \bar{q}_2\bar{q}_2' n)
\end{align}
and
\begin{align}
\label{b}
\mathfrak{B}(n;q_1,q_2,q_2')=\sum_{\alpha_2\bmod{q_2q_2'}}K_{q_2}(\alpha_2\ell\bar{q}_1^3)K_{q_2'}(-\alpha_2\ell\bar{q}_1^3)e_{q_2q_2'}(\alpha_2\bar{q}_1 n).
\end{align}\\

In the rest of this section we will analyse \eqref{b}. For $q_2=q_2'$ we use the bounds for the hyper-Kloosterman sums without trying to get extra cancellation in the some over $\alpha_2$. So we have
\begin{align*}
\mathfrak{B}(n;q_1,q_2,q_2')\ll Q_2^4.
\end{align*}
Now consider the case $q_2\neq q_2'$ (i.e. $(q_2,q_2')=1$). Then the sum further splits as
$$
\sum_{\alpha_2\bmod{q_2}}K_{q_2}(\alpha_2\ell\bar{q}_1^3)e_{q_2}(\alpha_2\bar{q}_1\bar{q}_2' n)\sum_{\alpha_2'\bmod{q_2'}}K_{q_2'}(-\alpha_2\ell\bar{q}_1^3)e_{q_2'}(\alpha_2\bar{q}_1\bar{q}_2 n).
$$
Consider the first sum. Opening the hyper-Kloosterman sum we get
$$
\sideset{}{^\star}\sum_{a,b\bmod{q_2}}e_{q_2}(a+b)\sum_{\alpha_2\bmod{q_2}}e_{q_2}(\alpha_2\ell\bar{q}_1^3\bar{a}\bar{b}+\alpha_2\bar{q}_1\bar{q}_2' n).
$$
Since $(\ell,q_2q_2')=1$ (as $\ell< Q_2$), the inner sum vanishes unless $(n,q_2)=1$. In particular this implies that if $q_2\neq q_2'$ then the zero frequency $n=0$ contribution vanishes. Now if $(n,q_2)=1$ then the above sum can be expressed in terms of the Kloosterman sum, namely
$$
q_2S(1,\ell \bar{q}_1^2q_2'\bar{n};q_2).
$$
We conclude the following lemma.

\begin{lemma}
\label{lem-char-1}
We have
\begin{align*}
\mathfrak{B}(n;q_1,q_2,q_2')\ll \begin{cases} Q_2^4 &\text{if}\;\;q_2=q_2'\\
Q_2^3 &\text{otherwise}.
\end{cases}
\end{align*}
Also $\mathfrak{B}(0;q_1,q_2,q_2')=0$ if $q_2\neq q_2'$.
\end{lemma}

\bigskip

%===============================================================================

\section{A character sum}
It remains to estimate the character sum \eqref{a}. %In the case $q_1|n$ we will be satisfied with the bound 
%$$
%\mathfrak{A}(n;q_1,q_2,q_2')\ll Q_1^3
%$$
%which follows from the Deligne bound for hyper-Kloosterman sums. Now 
We first consider the case where $q_1\nmid n$. Opening the hyper-Kloosterman sums we arrive at 
\begin{align*}
\mathfrak{A}(n;q_1,q_2,q_2')=\sideset{}{^\star}\sum_{a,b\bmod{q_1}}&\;\sideset{}{^\star}\sum_{c,d\bmod{q_1}}e_{q_1}(a+b+c+d)\\
&\times\sum_{\alpha\bmod{q_1}}e_{q_1}(\alpha_1\ell\bar{q}_2^3\bar{a}\bar{b}-\alpha_1\ell\bar{q_2'}^3\bar{c}\bar{d}+\alpha_1 \bar{q}_2\bar{q}_2' n).
\end{align*}
Executing the sum over $\alpha_1$ we obtain
\begin{align}
\label{mid-way}
q_1\mathop{\sideset{}{^\star}\sum_{a,b\bmod{q_1}}\;\sideset{}{^\star}\sum_{c,d\bmod{q_1}}}_{\ell\bar{q}_2^3\bar{a}\bar{b}-\ell\bar{q_2'}^3\bar{c}\bar{d}+ \bar{q}_2\bar{q}_2' n\equiv 0\bmod{q_1}}e_{q_1}(a+b+c+d).
\end{align}
From the congruence condition we can uniquely solve for $\bar{a}$. We get
\begin{align*}
\bar{a}\equiv \bar{\ell}q_2^3b(\ell\bar{q_2'}^3\bar{c}\bar{d}- \bar{q}_2\bar{q}_2' n) \bmod{q_1}.
\end{align*}
It follows that $cdnq_2'^2-\ell q_2$ is necessarily coprime with $q_1$. It follows that 
\begin{align*}
\mathfrak{A}(n;q_1,q_2,q_2')=q_1\sideset{}{^\star}\sum_{b\bmod{q_1}}\;\sideset{}{^\dagger}\sum_{c,d\bmod{q_1}}e_{q_1}(-\ell \bar{b}cd\bar{q}_2^2q_2'^3\overline{(cdnq_2'^2-\ell q_2)}+b+c+d),
\end{align*}
where the $\dagger$ implies that we have the restriction $(cd(cdnq_2'^2-\ell q_2),q_1)=1$. Let 
$$
\xi=cdnq_2'^2-\ell q_2.
$$
Note that since we are in the case $(q_1,n)=1$, $c$ can be uniquely determined from $(\xi,d)$. Moreover we have $(\xi(\xi+\ell q_2),q_1)=1$.\\

Thus we obtain that
\begin{align*}
\mathfrak{A}(n;q_1,q_2,q_2')=q_1\sideset{}{^\star}\sum_{b\bmod{q_1}}\;\sideset{}{^\star}\sum_{\substack{\xi,d\bmod{q_1}\\(\xi+\ell q_2,q_1)=1}}e_{q_1}(-\ell \bar{b}(\xi+\ell q_2)\bar{q}_2^2q_2'\bar{n}\overline{\xi}+b+(\xi+\ell q_2)\bar{q}_2'^2\bar{n}\bar{d}+d).
\end{align*}
Adding and subtracting the missing value $\xi=-\ell q_2$ to the sum it follows that
\begin{align*}
\mathfrak{A}(n;q_1,q_2,q_2')=q_1\sideset{}{^\star}\sum_{b\bmod{q_1}}\;\sideset{}{^\star}\sum_{\xi,d\bmod{q_1}}e_{q_1}(-\ell \bar{b}(\xi+\ell q_2)\bar{q}_2^2q_2'\bar{n}\overline{\xi}+b+(\xi+\ell q_2)\bar{q}_2'^2\bar{n}\bar{d}+d)+q_1.
\end{align*}
Let $\mathbf{x}=(x_1,x_2,x_3)$ and define the Laurent polynomial
$$
f(\mathbf{x})=-\ell \bar{q}_2^2q_2'\bar{n}\frac{1}{x_1}-\ell^2 \bar{q}_2q_2'\bar{n}\frac{1}{x_1x_2}+x_1+\bar{q}_2'^2\bar{n}\frac{ x_2}{x_3}+\ell q_2\bar{q}_2'^2\bar{n}\frac{1}{x_3}+x_3.
$$
We can now write
\begin{align*}
\mathfrak{A}(n;q_1,q_2,q_2')=q_1\sum_{\mathbf{x}\in (\mathbb{F}_{q_1}^\times)^3}e_{q_1}(f(\mathbf{x}))+q_1.
\end{align*}
Since none of the coefficients vanish modulo $q_1$, the Newton polyhedron $\Delta_\infty(f)$ of $f$ at infinity is given by the convex hull of the vectors $\mathbf{0}$, $\pm \mathbf{e}_1$, $\pm \mathbf{e}_3$, $-\mathbf{e}_1-\mathbf{e}_2$ and $\mathbf{e}_2-\mathbf{e}_3$. (Here $\mathbf{e}_i$ denotes the vector in $\mathbb{R}^3$ with $1$ at the $i$-th place and zeros elsewhere.) Clearly $\Delta_\infty(f)$ is of full dimension. Now we should check the non-degeneracy of $f$ with respect to $\Delta_\infty(f)$ (see \cite{DL}), i.e. we have to check that for any face $\tau$ of $\Delta_\infty$ not containing the origin, with associated Laurent polynomial $f_\tau$, the variety 
\begin{align}
\label{variety}
\frac{\partial}{\partial x_1}f_\tau =\frac{\partial}{\partial x_2}f_\tau=\frac{\partial}{\partial x_3}f_\tau=0
\end{align}
is empty. (Recall that for any $\tau$ the Laurent polynomial $f_\tau$ is defined by picking those terms from $f$ which have their corresponding index vector in $\tau$.) A moments reflection shows that we need to worry only about the case where $\tau$ has dimension $2$. In fact for \eqref{variety} to hold $\tau$ need to have at least four vertices. But any $\tau$ containing four vertices either contains the origin or has dimension three. Hence $f$ is non-degenerate. Then the main result of \cite{DL} implies that
\begin{align*}
\mathfrak{A}(n;q_1,q_2,q_2')\ll Q_1^{5/2}
\end{align*}
if $q_1\nmid n$.\\

Now consider the case where $q_1|n$. Then from \eqref{mid-way} we get
\begin{align*}
\mathfrak{A}(n;q_1,q_2,q_2')=q_1\sideset{}{^\star}\sum_{b\bmod{q_1}}\;\sideset{}{^\star}\sum_{c,d\bmod{q_1}}e_{q_1}(q_2'^3\bar{q}_2^3\bar{b}cd+b+c+d).
\end{align*}
The sum over $c$ now gives a Ramanujan sum. It follows that
\begin{align*}
\mathfrak{A}(n;q_1,q_2,q_2')=q_1^2\sideset{}{^\star}\sum_{b\bmod{q_1}}\;e_{q_1}(b-\bar{q}_2'^3q_2^3 b)-q_1.
\end{align*}
Consequently
\begin{align*}
\mathfrak{A}(n;q_1,q_2,q_2')\ll Q_1^2(q_1,q_2-q_2').
\end{align*}\\

\begin{lemma}
We have
\begin{align*}
\mathfrak{A}(n;q_1,q_2,q_2')\ll\begin{cases}Q_1^2(q_1,q_2-q_2') &\text{if}\;\;q_1|n;\\
Q_1^{5/2} &\text{otherwise}.
\end{cases}
\end{align*}
\end{lemma}

Combining with Lemma~\ref{lem-char-1}, we conclude the following.
\begin{corollary}
\label{cor}
We have
\begin{align*}
\mathfrak{C}(n;q_1,q_2,q_2')\ll\begin{cases}0 &\text{if}\;\;n=0\;\;\text{and}\;\;q_2\neq q_2';\\
Q^3Q_2 &\text{if}\;\;n=0\;\;\text{and}\;\;q_2=q_2';\\
Q^{5/2}\sqrt{Q_2}(q_1,n) &\text{if}\;\;n\neq 0\;\;\text{and}\;\;q_2\neq q_2';\\
Q^{5/2}Q_2^{3/2}(q_1,n) &\text{if}\;\;n\neq 0\;\;\text{and}\;\;q_2=q_2'.
\end{cases}
\end{align*}
\end{corollary}

\bigskip

%===================================================================================
\section{Conclusion}

Substituting the bounds from Corollary~\ref{cor} into \eqref{to-replace} we get
\begin{align*}
\mathfrak{E}_{q_1,q_2,q_2'}\ll \begin{cases} \sqrt{QX}Q^3+Q^{5/2+\varepsilon}Q_2^{3/2} &\text{if}\;\;q_2=q_2';\\
Q^{5/2+\varepsilon}\sqrt{Q_2} &\text{if}\;\;q_2\neq q_2'.
\end{cases}
\end{align*}
Substituting in \eqref{to-replace-1} we obtain
\begin{align*}
\mathfrak{E}(s)\ll Q^\varepsilon Q_1(\sqrt{QX}Q^3+Q^{5/2}Q_2^{3/2}).
\end{align*}
for $s=\varepsilon+it$ with $|t|\ll Q^\varepsilon$. Substituting this into Lemma~\ref{lem-s}, we get
$$
\mathfrak{S}\ll Q^{5/4+\varepsilon}(X^{1/4}Q^{1/2}Q_1^{1/2}+Q^{1/2}Q_2^{1/4}+(Q_1+Q_2)X^{1/2})
$$
Taking $X=Q^{-1/2+\delta}$ we get
$$
\mathfrak{S} \ll \left[Q^{13/8+\delta/4}Q_1^{1/2}+Q^{7/4}Q_2^{1/4}+(Q_1+Q_2)Q^{1+\delta/2}\right]Q^\varepsilon.
$$
Now we pick $Q_1=Q^{3/4-\delta}$. Then we have $Q_2=Q^{1/4+\delta}$, and the above expression yields
$$
\mathfrak{S} \ll \left[Q^{2-\delta/4}+Q^{7/4+1/16+2\delta}\right]Q^\varepsilon.
$$
(Note that we are not trying to obtain an optimal bound for the error term, as our main theorem does not require a strong asymptotic.) Picking $\delta=1/100$, and combining with Lemma~\ref{lem-f} (and \eqref{err-term-1}) we conclude Proposition~\ref{prop}.\\

%+++++++++++++++++++++++++++++++++++++++++++++++++++++++++++++++++++++++++++++++++++

\end{document}